\documentclass[11pt,a4paper]{article}
\usepackage[utf8]{inputenc}
\usepackage[english]{babel}
\usepackage{amssymb, amsmath, enumerate}
\usepackage{graphicx}
\usepackage{epstopdf}
\usepackage{subcaption}
\usepackage{authblk}
\usepackage{textcomp}

\newenvironment{definition}[1][Definition]{\begin{trivlist}
\item[\hskip \labelsep {\bfseries #1}]}{\end{trivlist}}

\newcommand{\qed}{\nobreak \ifvmode \relax \else
      \ifdim\lastskip<1.5em \hskip-\lastskip
      \hskip1.5em plus0em minus0.5em \fi \nobreak
      \vrule height0.75em width0.5em depth0.25em\fi}

\begin{document}
\title{Investment projects implementation with production facilities location 
 taking into account the environmental pollution}
\author[1]{O.A. Malafeyev\thanks{o.malafeev@spbu.ru}}
\author[1]{Y.E. Lakhina\thanks{juliala1401@gmail.com}}
\author[1]{N.D.Redinskikh \thanks{redinskich@yandex.ru}}
\affil[1]{Saint-Petersburg State University,  Russia}
\date{}
\maketitle
\begin{abstract}
	When implementing an investment projects a problem of production facilities location with taking into account the environmental pollution arises. A mathematical model of the industrial enterprises location in the region is formalized and studied in this paper. It is necessary to locate objects polluting the environment in such way as to maximize the total income of the players. The total income is calculated as income from the activities of enterprises, minus the funds that are spent on the recover damage to the environment. This problem is formalized as a non-cooperative game with  $n$ players that exploit common sources -- natural objects. Numerical example is solved.
\end{abstract}

\section{Introduction}

	Let us define a Cartesian coordinate system. The axis $OX$ is directed to the east, the axis $OY$ is directed strictly to the north.  In the region under consideration  there are natural objects $A_1,\ldots, A_m$. Each object is characterized by a pair of coordinates $(x_j, y_j)$.  Participants in the competitive process are objects: "pollutants" $B_1,\ldots,B_n$ and "contaminated"  $A_1,\ldots, A_m$. The number of  players is $n$, and  the number of natural objects is $m$.  
	
	The basics of the approaches and methodologies applied here can also be seen in [1-50].

\section{Main results}
The set of the player $i$ strategies :
\[
D_i=[u_i=(x_i,y_i), \underline{\rho} \leq \sqrt{(x_i-x_j)^2+(y_i-y_j)^2} \leq \overline{\rho}, \forall j \in m, j\neq i], i \in n
\]
\[
0\leq x_i \leq x_{max}, 0\leq y_i \leq y_{max}.
\]
where  $\underline{\rho} (\overline{\rho})$ -- is the minimum (maximum) distance between objects.

For each players strategy profile the players payoff functions are defined. The income of the player is equal to the value of its payoff function.

Let us take the income function for the $i$ player as follows:
\begin{eqnarray}\label{[FD]}
H_i(u)=\sum\limits_{j=1,j\neq i}^m\dfrac{L_{ij}}{\rho(i,j)}-\sum\limits_{j=1,j\neq i}^m\dfrac{Q_{ij}W_i}{2\pi\rho^2(i,j)}, i=\overline{1,n},
\end{eqnarray}
where  $n$  is the number of players, $m$ is the number of natural objects;
$\rho(i,j)=\sqrt{(x_i-x_j)^2+(y_i-y_j)^2}$  is the distance between the industrial enterprise $i$ and the natural object $j$.

The income from the activity of the enterprise $i$, received by the  player $i$ is:
\begin{eqnarray}\label{[D]}
\dfrac{L_{ij}}{\rho(i,j)},
\end{eqnarray}
where  $L_{ij} \geq 0$ is the amount of the loss, depending on the distance $\rho(i,j)$ between the objects $i$ and $j$. The greater the distance between the enterprise and the natural object, the lower of the income amount of the $i$ player will receive.

The amount of resources that the player $i$ spends on compensating for environmental damage to the natural object $j$:
\begin{eqnarray}\label{[Y]}
\dfrac{Q_{ij}W_i}{2\pi\rho^2(i,j)}.
\end{eqnarray}
$Q_{ij}$ is the weighting factor. It defines environmental damage that object $i$ causes   to object $j$, $W_i$ is the amount of harmful substances that the object $i$ emits into the environment.

If the object $i$ does not harm the object $j$, then $Q_{ij}=0$.

The function $H_i(u)$ in the domain of $(\underline{\rho} \leq \sqrt{(x_i-x_j)^2+(y_i-y_j)^2} \leq \overline{\rho})$ is continuous and has partial derivatives of the first and second order that are continuous in this region. Hence, it is smooth in the domain of the job.

\section{Numerical example}

Let's consider a numerical example.

A noncooperative game  $G=(N,\{X_i\}_{i\in N},\{H_i\}_{i\in N})$ is considered where  $N=3$ is the number of players, $X_i$ is the set of player $i$ strategies, $H_i$ is the player $i$ payoff function.
Let the number of natural objects be $m=5$. Let the region has an area of 15 square kilometers. We set  $\pi=3$. Let the natural objects $A_1, A_2, A_3, A_4, A_5$ in the region be arranged as follows:
\begin{table}[h!]
	\caption{\label{tab:canonsummary}Location of natural objects.}
	\begin{center}
		\begin{tabular}{|c|c|c|c|c|c|}
			\hline
			Natural objects & $A_1$ & $A_2$ & $A_3$ & $A_4$ & $A_5$ \\
			\hline
			Coordinate $x$ & $2$ & $5$ & $9$ & $14$ & $8$ \\
			Coordinate $y$ & $3$ & $9$ & $6$ & $1$ & $13$ \\
			\hline
		\end{tabular}
	\end{center}
\end{table}

For the enterprise of each player there is a permissible set of points where it is possible to build it, that is, each player $i$ has an acceptable set of strategies:
\begin{table}[h!]
	\caption{\label{tab:canonsummary}Admissible positions ($B_1, B_2, B_3$) for placing the company's player 1.}
	\begin{center}
		\begin{tabular}{|c|c|c|c|}
			\hline
			Industrial enterprise 1 & $B_1$ & $B_2$ & $B_3$ \\
			\hline
			Coordinate $x$ & $7$ & $1$ & $9$ \\
			Coordinate $y$ & $8$ & $2$ & $10$ \\
			\hline
		\end{tabular}
	\end{center}
\end{table}

\begin{table}[h!]
	\caption{\label{tab:canonsummary}Admissible positions ($C_1, C_2, C_3, C_4$) for placing the company's player 2.}
	\begin{center}
		\begin{tabular}{|c|c|c|c|c|}
			\hline
			Industrial enterprise 2 & $C_1$ & $C_2$ & $C_3$ & $C_4$\\
			\hline
			Coordinate $x$ & $6$ & $11$ & $5$ & $8$ \\
			Coordinate $y$ & $4$ & $15$ & $3$ & $15$ \\
			\hline
		\end{tabular}
	\end{center}
\end{table}

\begin{table}[h!]
	\caption{\label{tab:canonsummary}Admissible positions ($D_1, D_2$) for placing the company's player 3.}
	\begin{center}
		\begin{tabular}{|c|c|c|}
			\hline
			Industrial enterprise 3 & $D_1$ & $D_2$\\
			\hline
			Coordinate $x$ & $4$ & $6$  \\
			Coordinate $y$ & $12$ & $1$  \\
			\hline
		\end{tabular}
	\end{center}
\end{table}

\begin{table}[h!]
	\caption{\label{tab:canonsummary}The amount of loss $L_{ij}$ of player 1.}
	\begin{center}
		\begin{tabular}{|c|c|c|c|}
			\hline
			Industrial enterprise 1 & $B_1$ & $B_2$ & $B_3$\\
			\hline
			Natural object  $A_1$ & $10$ & $1$ & $13$ \\
			Natural object  $A_2$ & $4$ & $11$ & $8$ \\
			Natural object  $A_3$ & $5$ & $12$ & $6$ \\
			Natural object  $A_4$ & $13$ & $15$ & $14$ \\
			Natural object  $A_5$ & $9$ & $15$ & $6$ \\
			\hline
		\end{tabular}
	\end{center}
\end{table}

\begin{table}[h!]
	\caption{\label{tab:canonsummary}The amount of loss $L_{ij}$ of player  2.}
	\begin{center}
		\begin{tabular}{|c|c|c|c|c|}
			\hline
			Industrial enterprise 2 & $C_1$ & $C_2$ & $C_3$ & $C_4$\\
			\hline
			Natural object  $A_1$ & $5$ & $17$ & $2$ & $15$\\
			Natural object  $A_2$ & $7$ & $10$ & $8$ & $9$ \\
			Natural object  $A_3$ & $4$ & $13$ & $6$ & $12$ \\
			Natural object  $A_4$ & $11$ & $16$ & $13$ & $18$ \\
			Natural object  $A_5$ & $13$ & $3$ & $14$ & $1$ \\
			\hline
		\end{tabular}
	\end{center}
\end{table}

\begin{table}[h!]
	\caption{\label{tab:canonsummary}The amount of loss $L_{ij}$ of player 3.}
	\begin{center}
		\begin{tabular}{|c|c|c|}
			\hline
			Industrial enterprise 3 & $D_1$ & $D_2$ \\
			\hline
			Natural object  $A_1$ & $8$ & $3$ \\
			Natural object  $A_2$ & $1$ & $7$ \\
			Natural object  $A_3$ & $5$ & $4$ \\
			Natural object  $A_4$ & $10$ & $6$ \\
			Natural object  $A_5$ & $2$ & $9$ \\
			\hline
		\end{tabular}
	\end{center}
\end{table}

\begin{table}[h!]
	\caption{\label{tab:canonsummary}The weighting coefficient $Q_{ij}$, which determines the environmental damage caused by an industrial enterprise 1 to natural object $j$.}
	\begin{center}
		\begin{tabular}{|c|c|c|c|}
			\hline
			Industrial enterprise 1 & $B_1$ & $B_2$ & $B_3$\\
			\hline
			Natural object  $A_1$ & $1.15$ & $2.75$ & $1.45$ \\
			Natural object  $A_2$ & $1.5$ & $1.95$ & $2.15$ \\
			Natural object  $A_3$ & $1$ & $1.15$ & $1.05$ \\
			Natural object  $A_4$ & $2.2$ & $1.8$ & $2.9$ \\
			Natural object  $A_5$ & $1.9$ & $2.6$ & $1.4$ \\
			\hline
		\end{tabular}
	\end{center}
\end{table}

\begin{table}[h!]
	\caption{\label{tab:canonsummary}The weighting coefficient  $Q_{ij}$, which determines the environmental damage for player 2.}
	\begin{center}
		\begin{tabular}{|c|c|c|c|c|}
			\hline
			Industrial enterprise 2 & $C_1$ & $C_2$ & $C_3$ & $C_4$\\
			\hline
			Natural object  $A_1$ & $2.4$ & $1.96$ & $1.34$ & $2.05$ \\
			Natural object  $A_2$ & $1.67$ & $1.02$ & $1.73$ & $1.09$ \\
			Natural object  $A_3$ & $2.45$ & $1.75$ & $1$ & $2.05$ \\
			Natural object  $A_4$ & $1.85$ & $2.3$ & $1.6$ & $1.31$ \\
			Natural object  $A_5$ & $1.1$ & $2.7$ & $1.32$ & $1.09$ \\
			\hline
		\end{tabular}
	\end{center}
\end{table}

\begin{table}[h!]
	\caption{\label{tab:canonsummary} The weighting coefficient $Q_{ij}$, which determines the environmental damage for player 3.}
	\begin{center}
		\begin{tabular}{|c|c|c|}
			\hline
			Industrial enterprise 3 & $D_1$ & $D_2$ \\
			\hline
			Natural object  $A_1$ & $2.9$ & $1.25$ \\
			Natural object  $A_2$ & $1.05$ & $1.64$ \\
			Natural object  $A_3$ & $2.1$ & $1.36$ \\
			Natural object  $A_4$ & $1.9$ & $1.82$ \\
			Natural object  $A_5$ & $1.08$ & $1.6$ \\
			\hline
		\end{tabular}
	\end{center}
\end{table}

\begin{table}[h!]
	\caption{\label{tab:canonsummary}Amount of harmful substances  $W_i$, discarded enterprise $i$ to favorites}
	\begin{center}
		\begin{tabular}{|c|c|}
			\hline
			Industrial enterprise 1 & $60$ \\
			Industrial enterprise 2 & $15$ \\
			Industrial enterprise 3 & $35$ \\
			\hline
		\end{tabular}
	\end{center}
\end{table}

\begin{table}[h!]
	\caption{\label{tab:canonsummary}The influence of the choice of strategies of players 2 and 3 on the income of the player 1, when using 1 strategy.}
	\begin{center}
		\begin{tabular}{|c|c|c|c|c|}
			\hline
			Player positions 3/Player positions 2  & $C_1$ & $C_2$ & $C_3$ & $C_4$\\
			\hline
			$D_1$ & $0.17$ & $0.89$ & $0.25$ & $0.32$ \\
			$D_2$ & $1.01$ & $1.82$ & $1.54$ & $1.76$ \\
			\hline
		\end{tabular}
	\end{center}
\end{table}

\begin{table}[h!]
	\caption{\label{tab:canonsummary}The influence of the choice of strategies of players 2 and 3 on the income of the player 1, when using 2 strategy.}
	\begin{center}
		\begin{tabular}{|c|c|c|c|c|}
			\hline
			Player positions 3/Player positions 2  & $C_1$ & $C_2$ & $C_3$ & $C_4$\\
			\hline
			$D_1$ & $1.54$ & $1.82$ & $1.01$ & $1.76$ \\
			$D_2$ & $0.25$ & $0.89$ & $0.17$ & $0.32$ \\
			\hline
		\end{tabular}
	\end{center}
\end{table}

\begin{table}[h!]
	\caption{\label{tab:canonsummary}The influence of the choice of strategies of players 2 and 3 on the income of the player 1, when using 3 strategy.}
	\begin{center}
		\begin{tabular}{|c|c|c|c|c|}
			\hline
			Player positions 3/Player positions 2  & $C_1$ & $C_2$ & $C_3$ & $C_4$\\
			\hline
			$D_1$ & $0.32$ & $0.25$ & $0.89$ & $0.17$ \\
			$D_2$ & $1.76$ & $1.54$ & $1.82$ & $1.01$ \\
			\hline
		\end{tabular}
	\end{center}
\end{table}

\begin{table}[h!]
	\caption{\label{tab:canonsummary}The influence of the choice of strategies of players 1 and 3 on the income of the player 2, when using 1 strategy}
	\begin{center}
		\begin{tabular}{|c|c|c|c|}
			\hline
			Player positions 3/Player positions 1  & $B_1$ & $B_2$ & $B_3$\\
			\hline
			$D_1$ & $0.74$ & $0.85$ & $0.96$ \\
			$D_2$ & $0.32$ & $0.46$ & $0.57$ \\
			\hline
		\end{tabular}
	\end{center}
\end{table}

\begin{table}[h!]
	\caption{\label{tab:canonsummary}The influence of the choice of strategies of players 1 and 3 on the income of the player 2, when using 2 strategy.}
	\begin{center}
		\begin{tabular}{|c|c|c|c|}
			\hline
			Player positions 3/Player positions 1  & $B_1$ & $B_2$ & $B_3$\\
			\hline
			$D_1$ & $0.46$ & $0.57$ & $0.32$  \\
			$D_2$ & $0.85$ & $0.96$ & $0.74$  \\
			\hline
		\end{tabular}
	\end{center}
\end{table}

\begin{table}[h!]
	\caption{\label{tab:canonsummary}The influence of the choice of strategies of players 1 and 3 on the income of the player 2, when using 3 strategy.}
	\begin{center}
		\begin{tabular}{|c|c|c|c|}
			\hline
			Player positions 3/Player positions 1  & $B_1$ & $B_2$ & $B_3$\\
			\hline
			$D_1$ & $0.85$ & $0.74$ & $0.96$  \\
			$D_2$ & $0.46$ & $0.32$ & $0.57$  \\
			\hline
		\end{tabular}
	\end{center}
\end{table}

\begin{table}[h!]
	\caption{\label{tab:canonsummary}The influence of the choice of strategies of players 1 and 3 on the income of the player 2, when using 4 strategy.}
	\begin{center}
		\begin{tabular}{|c|c|c|c|}
			\hline
			Player positions 3/Player positions 1  & $B_1$ & $B_2$ & $B_3$\\
			\hline
			$D_1$ & $0.46$ & $0.57$ & $0.32$  \\
			$D_2$ & $0.85$ & $0.96$ & $0.74$  \\
			\hline
		\end{tabular}
	\end{center}
\end{table}

\begin{table}[h!]
	\caption{\label{tab:canonsummary}The influence of the choice of strategies of players 1 and 2 on the income of the player 3, when using 1 strategy.}
	\begin{center}
		\begin{tabular}{|c|c|c|c|c|}
			\hline
			Player positions 1/Player positions 2  & $C_1$ & $C_2$ & $C_3$ & $C_4$\\
			\hline
			$B_1$ & $0.65$ & $0.45$ & $0.74$ & $0.12$ \\
			$B_2$ & $1.63$ & $1.6$ & $1.7$ & $1.5$ \\
			$B_3$ & $1.25$ & $0.96$ & $1.32$ & $0.8$ \\
			\hline
		\end{tabular}
	\end{center}
\end{table}

\begin{table}[h!]
	\caption{\label{tab:canonsummary}The influence of the choice of strategies of players 1 and 2 on the income of the player 3, when using 2 strategy.}
	\begin{center}
		\begin{tabular}{|c|c|c|c|c|}
			\hline
			Player positions 1/Player positions 2  & $C_1$ & $C_2$ & $C_3$ & $C_4$\\
			\hline
			$B_1$ & $0.45$ & $0.74$ & $0.12$ & $0.65$ \\
			$B_2$ & $0.96$ & $1.32$ & $0.8$ & $1.25$ \\
			$B_3$ & $1.6$ & $1.7$ & $1.5$ & $1.63$ \\
			\hline
		\end{tabular}
	\end{center}
\end{table}

We calculate the value of the payoff function of each player according to the formula (\ref{[FD]}) taking into account the data specified above. We get the following payoff matrices for players 1, 2 and 3, where the lines are player's strategy 1, the columns are the strategies of player 2, the matrices are the strategies of player 3:
\[\begin{bmatrix}
(0.444; 3.931; 1.007) & (2.326; 2.565; 0.186) & (0.654; 4.220; 2.633) & (0.836; 3.759; 1.487) \\
(5.339; 4.515; 0.697) & (6.309; 3.178; 2.525) & (3.501; 3.674; 2.323) & (6.101; 4.658; 2.044) \\
(1.154; 5.100; 1.146) & (0.902; 1.784; 2.478) & (3.210; 4.766; 1.936) & (0.613; 2.615; 1.239) \\
\end{bmatrix}\]

\[\begin{bmatrix}
(2.640; 1.700; 1.201) & (4.757; 4.739; 1.735) & (4.025; 2.284; 2.135) & (4.600; 6.946; 4.537) \\
(0.867; 2.444; 1.975) & (3.085; 5.352; 2.562) & (0.589; 1.589; 3.336) & (1.109; 7.845; 4.003) \\
(6.348; 3.028; 0.320) & (5.554; 4.126; 3.523) & (6.564; 2.830; 4.270) & (3.643; 6.047; 4.350) \\
\end{bmatrix}\]

\subsection{The Nash Equilibrium}

In the book \cite{Zub} it is presented the construction of the Nash equilibrium strategy profile. Since in our problem each of the three players has a finite number of strategies, and each competitive strategy profile corresponds to the set of the income functions values of the players, then in this game there is at least one Nash equilibrium strategy profile in the mixed strategies. Let us find the equilibrium strategy profile in this example. Let us recall the definition of the equilibrium strategy profile. 
\begin{definition}
A set of the strategies $u^*=(u_1^*,\dots,u_k^*), (u_i^*)\in D_i$ is called a Nash equilibrium strategy profile if for any strategy  $u_i\in D_i$ the following inequality is valid
\begin{eqnarray}\label{[7]}
H_i(u^*)\geq H_i(u^*||u_i), i=\overline{1,k},
\end{eqnarray}
where $(u^*||u_i)=(u_1^*,\dots,u_{i-1}^*,u_{i+1}^*,\dots,u_k^*)$.
\end{definition}

From the definition of equilibrium it follows that an agent $i$ one-sided deviation can only lead  to a decrease in his income.

In this problem, the Nash equilibrium strategy profile is searched as follows: in the first step let us fix the strategy 1 of player 2 (column 1) and strategy 1 of player 3 (matrix 1). Then let us go through all the values (three) of the payoff function of player 1 and choose the largest of them.

In the second step let us fix strategy 2 of player 2 (column 2) and strategy 1 of player 3 (matrix 1). Similarly, let us find the largest value of the payoff function of player 1.

Thus, fixing the strategies of 2 and 3 players, let us find the greatest value of the payoff function of the player 1 in each column of each matrix.

After that let us search the largest value of the payoff function of player 2. Let us fix player's strategy 3 (matrix) and player's strategy 1 (fixed-matrix string). Further let us compare the values player's payoff function 2 (on matrix columns). Let us find the largest value in each row of both matrices. 

Similarly, by fixing the strategies of player 1 (row) and player 2 (column) let us choose the largest value of the payoff function of player 3 by its two strategies (matrices).

Finally, let us examine the intersection of the chosen win values for each player. This is the Nash equilibrium strategy profile.

Following this algorithm, let us obtain the following Nash equilibrium point: (4.600; 6.946; 4.537).

\subsection{A compromise strategy profile}

Let us find a compromise strategy profile in this problem.
The set of compromise strategies profile $M_i=\max\{H_i(u),u\in D\}$ is defined as follows:
$C_k=\{u\in D| \max \limits_i(M_i-H_i(u^{'}))\leq\max \limits_i(M_i-H_i(u)), \forall u\in D\}$.

In other words, a compromise strategy profile is a strategy profile in which the largest deviation of the payoff function of one of the players from its maximum value for $i$ is not greater than the largest deviation of the $i$ payoff function to one of the players from the maximum value in any other strategy profile.

The compromise strategy profile is searched as follows in this problem: for each player let us find the maximum value of its payoff function $M_i$:
\begin{equation}
M_1=6.564; M_2=7.845; M_3=4.537.
\end{equation}
Thus, the ideal vector is $(M_1,M_2,M_3)=(6.564,7.845,4.537)$.

Let us calculate the maximum values of the residuals for each strategy profile $x\in X =B*C*D$. Here $|X|=3*4*2=24$. Now let us calculate $x\in X$ maximum residual value $\delta (x)$:
\begin{equation}
(B_1,C_1,D_1): \delta (1)=M_1-H_1(1)=6.12;
\end{equation}
\begin{equation}
(B_1,C_2,D_1): \delta (2)=M_2-H_2(2)=5.28;
\end{equation}
\begin{equation}
(B_1,C_3,D_1): \delta (3)=M_1-H_1(1)=5.91;
\end{equation}
\begin{equation}
(B_1,C_4,D_1): \delta (4)=M_1-H_1(1)=5.728;
\end{equation}
\begin{equation}
(B_2,C_1,D_1): \delta (5)=M_3-H_3(1)=3.837;
\end{equation}
\begin{equation}
(B_2,C_2,D_1): \delta (6)=M_2-H_2(2)=4.667;
\end{equation}
\begin{equation}
(B_2,C_3,D_1): \delta (7)=M_2-H_2(3)=4.171;
\end{equation}
\begin{equation}
(B_2,C_4,D_1): \delta (8)=M_2-H_2(4)=3.187;
\end{equation}
\begin{equation}
(B_3,C_1,D_1): \delta (9)=M_1-H_1(3)=5.41;
\end{equation}
\begin{equation}
(B_3,C_2,D_1): \delta (10)=M_2-H_2(2)=6.061;
\end{equation}
\begin{equation}
(B_3,C_3,D_1): \delta (11)=M_1-H_1(3)=3.354;
\end{equation}
\begin{equation}
(B_3,C_4,D_1): \delta (12)=M_1-H_1(3)=5.951;
\end{equation}
\begin{equation}
(B_1,C_1,D_2): \delta (13)=M_2-H_2(1)=6.145;
\end{equation}
\begin{equation}
(B_1,C_2,D_2): \delta (14)=M_2-H_2(2)=3.106;
\end{equation}
\begin{equation}
(B_1,C_3,D_2): \delta (15)=M_2-H_2(3)=5.561;
\end{equation}
\begin{equation}
(B_1,C_4,D_2): \delta (16)=M_1-H_1(1)=1.964;
\end{equation}
\begin{equation}
(B_2,C_1,D_2): \delta (17)=M_1-H_1(2)=5.697;
\end{equation}
\begin{equation}
(B_2,C_2,D_2): \delta (18)=M_1-H_1(2)=3.479;
\end{equation}
\begin{equation}
(B_2,C_3,D_2): \delta (19)=M_2-H_2(3)=6.256;
\end{equation}
\begin{equation}
(B_2,C_4,D_2): \delta (20)=M_1-H_1(2)=5.455;
\end{equation}
\begin{equation}
(B_3,C_1,D_2): \delta (21)=M_2-H_2(1)=4.817;
\end{equation}
\begin{equation}
(B_3,C_2,D_2): \delta (22)=M_2-H_2(2)=3.719;
\end{equation}
\begin{equation}
(B_3,C_3,D_2): \delta (23)=M_2-H_2(3)=5.015;
\end{equation}
\begin{equation}
(B_3,C_4,D_2): \delta (24)=M_1-H_1(3)=2.921;
\end{equation}
Among them, let us find the minimum value (maximal discrepancy) for all strategies profile $x\in X$:
\begin{equation}
(B_1,C_4,D_2): \delta (16)=M_1-H_1(1)=1.964,
\end{equation}
that is, the desired compromise strategy profile is the strategy profile $(B_1,C_4,D_2)$.
The values of the players winning features in this strategy profile are the following $(4.600; 6.946;4.537)$.

\section{Conclusion}

The two principles of optimality, Nash equilibrium and compromise solution, are considered in this paper. The compromise solution in the model is the strategy profile $(B_1,C_4,D_2)$.

\section{Acknowledgements}
The work is partly supported by work RFBR No. 18-01-00796.

\end{document}